\documentclass[11pt,titlepage]{article}

\usepackage{bbm}
\usepackage{rotating, rotfloat}
\usepackage[round]{natbib} 
\usepackage{graphics, graphicx, hyperref, url, amsmath, amsthm, amsfonts, enumerate, epstopdf, booktabs,verbatim, subfigure}
\hypersetup{
    pdfstartview={FitH},    
    colorlinks=true,        
    linkcolor=blue,         
    citecolor=blue,         
    filecolor=blue,         
    urlcolor=blue           
}

     \def\etal{{\itshape et al.}}
\def\beqn{\begin{eqnarray}} \def\eeqn{\end{eqnarray}} 

\def\u{\textbf{\emph{u}}} \def\x{\textbf{\emph{x}}}  \def\z{\textbf{\emph{z}}}  \def\vv{\textbf{\emph{v}}}
\def\y{\textbf{\emph{y}}}

\newcommand{\HM}{\mbox{\it HM}}

\newcommand{\bec}{\begin{center}}
\newcommand{\enc}{\end{center}}
\newcommand{\bee}{\begin{eqnarray*}}
\newcommand{\ene}{\end{eqnarray*}}
\newcommand{\beq}{\begin{equation}}
\newcommand{\eeq}{\end{equation}}

\begin{document}

\title{\bf Finite sample breakdown point of Tukey's halfspace median}

\author {{Xiaohui Liu$^{a, c, d}$,\footnote{Corresponding author's email: csuliuxh912@gmail.com.}
        Yijun Zuo$^b$, Qihua Wang$^{c, e}$,}\\ \\
        {\em\footnotesize $^a$ School of Statistics, Jiangxi University of Finance and Economics, Nanchang, Jiangxi 330013, China}\\
         {\em\footnotesize $^b$ Department of Statistics and Probability, Michigan State University, East Lansing, MI, 48823, USA}\\
         {\em\footnotesize $^c$ Academy of Mathematics and Systems Sciences, Chinese Academy of Sciences, Beijing 100190, China}\\
         {\em\footnotesize $^d$ Research Center of Applied Statistics, Jiangxi University of Finance and Economics, Nanchang,}\\ {\em\footnotesize Jiangxi 330013, China}\\
         {\em\footnotesize $^e$ Institute of Statistical Science, Shenzhen University, Shenzhen 518006, China}\\
}

\maketitle

\begin{center}
{\sc Summary}
\end{center}

Tukey's halfspace median ($\HM$), servicing as the {multivariate} counterpart of the univariate median, has been introduced and extensively studied in the literature. It is supposed and expected to preserve robustness property (the most outstanding property) of the univariate median. One of prevalent quantitative assessments of robustness is finite sample breakdown point (FSBP). Indeed, the FSBP of many multivariate medians have been identified, except for the most prevailing one---the Tukey's halfspace median. This paper presents a precise result on FSBP for Tukey's halfspace median. The result here depicts the complete prospect of the global robustness of $\HM$ in the \emph{finite sample} practical scenario, revealing the dimension effect on the breakdown point robustness  and complimenting the existing \emph{asymptotic} breakdown point result.
\vspace{2mm}

\vspace{2mm}


\setlength{\baselineskip}{1.5\baselineskip}

\vskip 0.1 in
\section{Introduction}
\paragraph{}
\vskip 0.1 in \label{Introduction}

Robustness (as an insurance) is one of the most desirable properties for any statistical procedures.
The most outstanding feature of univariate median is its robustness.  Indeed, among all \emph{translation equivariant}  location estimators, it has the best possible \emph{breakdown point} \citep{Don1982} (and minimum \emph{maximum bias} if underlying distribution has unimodal symmetric density \citep{Hub1964}).
\medskip

It is very much desirable to extend the univariate median to multidimensional settings and meanwhile inherit/preserve its outstanding robustness for multidimensional data. In fact, the earliest attempt of this type of extension was made at least one century ago \citep{Web1909}. Oja's median \citep{Oja1983} is another promising extension.
\medskip

On the other hand, defining the multi-dimensional median as the deepest point of the underlying multidimensional data is an obvious natural approach. Serving this purpose, general notions of data depth have been proposed and studied \citep{ZS2000}. The main goal of data depth is to provide a center-outward ordering of multidimensional observations. Multivariate medians as the deepest point (the generalization of the univariate median) therefore have been naturally introduced and examined. Among the depth induced multidimensional medians, Tukey's halfspace median \citep{Tuk1975} is the most prominent and prevailing.  Robustness is of course the main targeted property to be shared by all depth induced medians.
\medskip

There are many ways to measure the robustness of a statistical procedure (especially location estimators). Among others, \emph{maximum bias}, \emph{influence function} and \emph{finite sample breakdown point} (FSBP) are the most standard gauges. FSBP by far is the most prevailing quantitative assessment of robustness due to its plain definition (without involvement of probability/randomness concept).
\medskip

The concept of breakdown point was introduced by \cite{Hod1967} and Hampel (Ph. D. dissertation (1968), Univ. California, Berkeley) and extended by \cite{Ham1971} and developed further by, among others, \cite{Hub1981}. It has proved most successful in the context of location, scale and regression problems. Finite sample version of breakdown point has been proposed, promoted and popularized by \cite{Don1982} and \cite{DH1983} (DH83).
\medskip

The seminar paper of \cite{DG1992} (hereafter DG92) was devoted to extensively study the FSBP of multivariate location estimators including Tukey's halfspace depth induced location estimators, especially the halfspace median ($\HM$).  Specifically, DG92 established FSBP for many location estimators, including the lower bound of FSBP for halfspace median.  However, lower bound contains much scarce information about FSBP of $\HM$.  What is the exact FSBP of $\HM$ is still an open question.
\medskip

\cite{AY2002} and \cite{CT2002} have pioneered in studying the maximum bias of $\HM$. It is found that the asymptotic breakdown point of $\HM$ is $1/3$, as also given in DG92 (see also \cite{Che1995} and \cite{Miz2002}). The latter result however is obtained restricted to a sub-class of distributions (the absolutely continuous centrosymemetric ones) in the maximum bias definition, and when sample size $n$ approaches to the infinity.
Ironically, DG92 only provided the \emph{asymptotic} breakdown point for $\HM$ and uncharacteristically left its FSBP open.
Furthermore, the former does not provide any clue of the dimensional effect on the breakdown robustness and its behavior in finite sample practical scenario. To address this issue and provide a definite answer is the main objective of this manuscript.
\medskip

Let's end this section with some definitions.
A location statistical functional $T$ in $\mathcal{R}^d$ ($d\ge 1$) is said to be \emph{affine equivariant} if
\begin{eqnarray*}
 T(\Sigma \mathcal{X}^n + \textbf{\emph{b}}) = \Sigma T(\mathcal{X}^n) + \textbf{\emph{b}}, \label{AE}
\end{eqnarray*}
for any $d \times d$ nonsingular matrix $\Sigma$ and $\textbf{\emph{b}} \in R^d$, where $\Sigma \mathcal{X}^n + \textbf{\emph{b}} =\{\Sigma X_1 + \textbf{\emph{b}}, \cdots, \Sigma X_n + \textbf{\emph{b}}\}$ and $X_1,\cdots, X_n$ is a given random sample  in $\mathcal{R}^d$ (denote $\mathcal{X}^n = \{X_1,\cdots, X_n\}$ hereafter). When $d=1$ and $\Sigma = 1$, we call $T$ is \emph{translation equivariant}.

\vskip 3mm

Define the depth of a point $\x$ with respect to $\mathcal{X}^n \subset \mathcal{R}^d$ 
as
\begin{eqnarray*}
    D(\x, \mathcal{X}^n) = \inf_{\u \in \mathcal{S}^{d-1}} P_n(\u^\top X \leq \u^\top \x),
\end{eqnarray*}
where $\mathcal{S}^{d-1} = \{\z \in \mathcal{R}^d: \|\z\| = 1\}$, and $P_n$ denotes the empirical probability measure
and $\|\cdot\|$ stands for Euclidean norm. Denote
\begin{eqnarray*}
    \mathcal{M}(\mathcal{X}^n) = \{\x: D(\x, \mathcal{ X}^n) = \lambda^*(\mathcal{X}^n)\},
\end{eqnarray*}
where $\lambda^*(\mathcal{X}^n) = \sup_\x D(\x, \mathcal{X}^n)$.
Tukey's halfspace median is defined as, 
\begin{eqnarray*}
    T^*(\mathcal{X}^n) = \textbf{Ave}\left\{\x \in \mathcal{M}(\mathcal{X}^n) \right\},
\end{eqnarray*}
i.e., the average of all points that maximize $D(\x, \mathcal{X}^n)$. Clearly, when $d = 1$, $T^*(\mathcal{X}^n)$ reduces to the univariate sample median.
\medskip

The finite sample \emph{additional breakdown point} (ABP) of a location estimator $T$ at the given sample
$\mathcal{X}^n$ is defined as
\begin{eqnarray*}
ABP(T,\mathcal{X}^n) = \min_{1\le m\le n}\left\{\frac{m}{n+m}: \sup_{\mathcal{Y}^{m}}\left\|T(\mathcal{X}^n\cup \mathcal{Y}^{m})- T(\mathcal{X}^n)\right\| =\infty\right\},
\end{eqnarray*}
where $\mathcal{Y}^{m}$
denotes an arbitrary contaminating sample of size $m$, adjoining to $\mathcal{X}^n \subset \mathcal{R}^d$.
Namely, the
ABP of an estimator is the minimum additional fraction which could drive
the estimator beyond any bound. It is readily seen that the ABP of the sample
mean and the univariate median are $1/(n+1)$ and $1/2$, respectively. The latter is the best that one can expect for any translation equivariant location estimator \citep{Don1982}.
\vskip 3mm

\emph{Additional} breakdown point is one of the forms of the finite sample breakdown point notion, \emph{replacement} breakdown point (RBP) is the other one (DH83), where instead of adding contaminating points to $\mathcal{X}^n$, replacing $m$ points of $\mathcal{X}^n$ by $m$ arbitrary points. Some prefer RBP since it is arguably more close to the contamination in reality.  The two are actually equivalent in the sense of \cite{Zuo2001}. Further discussions on FSBP concept could be found in DH83 and \cite{LR1991}. By FSBP we mean ABP in the sequel.
\medskip

We anticipate that the approach and results here may be extended for the investigating the FSBP of estimators that are related to Tukey's halfspace depth function such as multiple-output quantile regression estimators \citep{HPS2010}, the maximum regression depth estimator \citep{RH1999}, and probably the functional halfspace depth estimators \citep{LR2009}.
\medskip

The remainder of this article is organized as follows. Section~\ref{lemmas} presents three preliminary lemmas for the main results, which will be proved in Section~\ref{main}. The article ends in Section~\ref{Conclusion} with some concluding remarks.

\vskip 0.1 in
\section{Three preliminary lemmas}\label{lemmas}
\paragraph{}
\vskip 0.1 in

Since the proof of the main result is rather complicated and long, we divide it into several parts and present some of them as lemmas. In this section, three preliminary lemmas are established. They play important roles in the proof of the main results.
\vskip 3mm

Without loss of generality, we assume that $\mathcal{X}^n$ is in general position (\textbf{IGP} hereafter) throughout this paper. That is, no more than $d$ sample points lie in a $(d-1)$-dimensional hyperplane. This assumption is common in the literature involving statistical depth functions and breakdown point robustness \citep{DG1992, MLB2009}. Since $\HM$ reduces to the sample median and its FSBP is known as $1/2$ for $d = 1$, we only focus on $d \ge 2$ in the sequel.
\medskip

For $\mathcal{X}^n$, under the \textbf{IGP} assumption, there must exist $N_n^d = {n \choose d}$ unit vectors, say $\mu_j \in \mathcal{S}^{d-1}$, $j = 1, 2, \cdots, N_n^d$, such that they are respectively normal to $N_n^d$ hyperplanes with each of which passing through $d$ observations. Since $N_n^d$ is finite when $n$ and $d$ are fixed, we can find a unit vector, say $\u$, satisfying
\begin{eqnarray}
\label{eqnu0}
    \u^\top \mu_j \neq 0, \quad \text{ for } \forall j \in \left\{ 1, 2, \cdots, N_n^d \right\}.
\end{eqnarray}

For simplicity, for any given $\u \in \mathcal{S}^{d-1}$, in the sequel we denote
\begin{eqnarray}
\label{mathbbA}
    \mathbb{A}_{\u} = (\u_1, \u_2, \cdots, \u_{d-1}),
\end{eqnarray}
where $\u_1, \u_2, \cdots, \u_{d-1}$ are orthogonal to $\u$, and together with $\u$, they form a set of standard basis vectors of $\mathcal{R}^d$. \emph{Remarkably}, although the choice of $\u_1, \u_2, \cdots, \u_{d-1}$ is not unique when $d \ge 2$, this fact does not effect the proofs presented in the rest of this paper due to the affine equivariance of \emph{HM} and its related Tukey depth, nevertheless. Hence, we pretend that $\mathbb{A}_{\u}$ is unique in the sequel.
\medskip

Write $\mathbf{x}_i = \mathbb{A}_\u^\top X_i$, $i = 1, 2, \cdots, n$, and $\mathbf{X}_\u^n = \{\mathbf{x}_1, \mathbf{x}_2, \cdots, \mathbf{x}_n\}$, call it $\mathbb{A}_{\u}$-\emph{projections of $\mathcal{X}^n$} hereafter. It will greatly facilitate our discussion, if $\mathbf{X}_{\u}^n$ is still in general position. Fortunately, Lemma 1 provides a positive answer.
\medskip

\textbf{Lemma 1}.  Suppose $\mathcal{X}^n$ is \textbf{IGP}, and $\u \in \mathcal{S}^{d-1}$ satisfies display \eqref{eqnu0}. Then $\mathbf{X}_\u^n = \{\mathbf{x}_1, \mathbf{x}_2$, $\cdots, \mathbf{x}_n\} \subset \mathcal{R}^{d-1}$ is in general position too when $d \ge 2$.
\medskip

\textbf{Proof}. If $\mathbf{X}_\u^n$ is not \textbf{IGP}, then there must exist a $(d-2)$-dimensional hyperplane $\mathcal{P}_1$ containing at least $d$ $\mathbb{A}_{\u}$-projections, say $\mathbf{x}_{i_1}, \mathbf{x}_{i_2}, \cdots, \mathbf{x}_{i_k}$ ($k\ge d$). Let $\mathbf{v} \in \mathcal{S}^{d-2}$ be the normal vector of $\mathcal{P}_1$. Then, we have
\begin{eqnarray*}
    \mathbf{v}^\top \mathbf{x}_{i_1} = \mathbf{v}^\top \mathbf{x}_{i_2} = \cdots = \mathbf{v}^\top \mathbf{x}_{i_k}.
\end{eqnarray*}
Recalling the definition of $\mathbf{x}_i$, 
 we further obtain
\begin{eqnarray}
\label{eqn1}
    (\mathbb{A}_\u \mathbf{v})^\top X_{i_1} = (\mathbb{A}_\u \mathbf{v})^\top X_{i_2} = \cdots = (\mathbb{A}_\u \mathbf{v})^\top X_{i_k}.
\end{eqnarray}

Write $\widetilde{\vv} = \mathbb{A}_\u \mathbf{v}$. Clearly, $\widetilde{\vv} \in \mathcal{R}^{d}$ and $(\widetilde{\vv})^\top \widetilde{\vv} = 1$, namely, $\widetilde{\vv} \in \mathcal{S}^{d-1}$. Hence, \eqref{eqn1} implies that one can find a $(d-1)$-dimensional hyperplane, with $\widetilde{\vv}$ being its normal vector, that passes through $k$ observations. This contradicts with the \textbf{IGP} assumption of $\mathcal{X}^n$ if $k > d$. When $k = d$, \eqref{eqn1} implies $\widetilde{\vv} \in \{\mu_j\}_{j=1}^{N_n^d}$. This obviously contradicts with the fact that $\u$ satisfies \eqref{eqnu0} due to $\widetilde{\vv}^\top \u = \mathbf{v}^\top \mathbb{A}_\u^\top \u = 0$. This completes the proof of this lemma. \hfill $\Box$
\medskip

\textbf{Remark 2.1} In fact, $\mathcal{X}^n$ is \textbf{IGP} if and only if $\mathbf{X}_\u^n$ is \textbf{IGP} for $\u$ satisfying \eqref{eqnu0}.
\medskip

To derive the FSBP of \emph{HM}, we need to investigate the maximum Tukey depth with respect to the $\mathbb{A}_\u$-projections of $\mathcal{X}^n$. The following Lemma 2 will play an important role during this process.
\medskip

We formally introduce some additional necessary notations as follows. For $\forall \x, \y\in \mathcal{R}^d$, let
\begin{eqnarray*}
    \mathcal{U}_{\x} &=& \{\u\in \mathcal{S}^{d-1}: ~P_n(\u^\top X \leq \u^\top \x) = D(\x, \mathcal{X}^n)\}
\end{eqnarray*}
be the set of all optimal vectors of $\x$ which realize the depth at $\x$ with respect to $\mathcal{X}^n$, and
\begin{eqnarray*}
    \mathcal{H}_{\x, \y} &=& \{\u\in \mathcal{S}^{d-1}: ~\u^\top \x < \u^\top \y\}
\end{eqnarray*}
the hemisphere determined by $\{\x,\y\}$. Furthermore, for $\forall \z \in \mathcal{M}(\mathcal{X}^n)$, let
\begin{eqnarray*}
    \mathcal{A}_{\z} &=& \{\x \in \mathcal{M}(\mathcal{X}^n): \mathcal{U}_\x \cap \mathcal{H}_{\x, \z} \neq \emptyset\},\\
    \mathcal{B}_{\z} &=& \{\x \in \mathcal{M}(\mathcal{X}^n): \mathcal{U}_\x \cap \mathcal{H}_{\x, \z} = \emptyset\},\\
    \widetilde{\mathcal{B}}_{\z} &=& \mathcal{B}_{\z} \setminus \{\z\}.
\end{eqnarray*}
Obviously, (i) $\mathcal{M}(\mathcal{X}^n) = \mathcal{A}_{\z} \cup \mathcal{B}_{\z}$, (ii) $\mathcal{A}_{\z} \cap \mathcal{B}_{\z} = \emptyset$, (iii) $\z \in \mathcal{B}_{\z}$ because $\mathcal{H}_{\z, \z} = \emptyset$. 
\medskip

For $\mathcal{U}_{\x}$, $\mathcal{B}_{\x}$ and $\widetilde{\mathcal{B}}_{\z}$, Lemma 2 below depicts several important properties of them.
\medskip

\textbf{Lemma 2}. Suppose $\mathcal{X}^n$ is \textbf{IGP} and $\mathcal{M}(\mathcal{X}^n)$ is of affine dimension $d$. For $\forall \z_1$ lying in the interior of $\mathcal{M}(\mathcal{X}^n)$, if $\widetilde{\mathcal{B}}_{\z_1} \neq \emptyset$, then:
\begin{enumerate}
  \item[] (\textbf{o1}) for $\forall \u \in \mathcal{U}_{\z_1}$ and $\forall \z \in \widetilde{\mathcal{B}}_{\z_1}$, we have $\u^\top \z \ge \u^\top \z_1$.

  \item[] (\textbf{o2}) for $\forall \z \in \mathcal{B}_{\z_1}$, we have $\mathcal{U}_\z \subset \mathcal{U}_{\z_1}$.

  \item[] (\textbf{o3}) for $\forall \z \in \widetilde{\mathcal{B}}_{\z_1}$, we have $\mathcal{B}_\z \subset \mathcal{B}_{\z_1}$, but $\z_1 \notin \mathcal{B}_\z$. That is, $\mathcal{B}_\z \subset \widetilde{\mathcal{B}}_{\z_1}$.
\end{enumerate}

\textbf{Proof}. (\textbf{o1}). If not, $\u^\top \z < \u^\top \z_1$. Then $\lambda^*(\mathcal{X}^n) \leq P_n (\u^\top X \leq \u^\top \z) \leq P_n(\u^\top X \leq \u^\top \z_1) = \lambda^*(\mathcal{X}^n)$, resulting in $\u \in \mathcal{U}_{\z} \cap \mathcal{H}_{\z, \z_1}$, and hence contradicting with $\z \in \mathcal{B}_{\z_1}$.
\medskip

(\textbf{o2}). By definition, $\z \in \mathcal{B}_{\z_1}$ implies $\mathcal{U}_{\z} \cap \mathcal{H}_{\z, \z_1} = \emptyset$. Hence, for $\forall \u \in \mathcal{U}_{\z}$, we have $\u^\top \z \ge \u^\top \z_1$. Then $\lambda^*(\mathcal{X}^n) \leq P_n (\u^\top X \leq \u^\top \z_1) \leq P_n(\u^\top X \leq \u^\top \z) = \lambda^*(\mathcal{X}^n)$. That is, $\u\in \mathcal{U}_{\z_1}$, and hence $\mathcal{U}_{\z} \subset \mathcal{U}_{\z_1}$.
\medskip

(\textbf{o3}). For $\forall \x \in \mathcal{B}_\z$, $\mathcal{U}_\x \bigcap \mathcal{H}_{\x, \z} = \emptyset$. Hence, $\u_\x^\top (\x - \z) \ge 0$ for $\forall \u_\x \in \mathcal{U}_\x$. Next, by (\textbf{o2}), we have $\mathcal{U}_\x \subset \mathcal{U}_{\z} \subset \mathcal{U}_{\z_1}$, which implies $\u_\x^\top \x \ge \u_\x^\top \z \ge \u_\x^\top \z_1$ by (\textbf{o1}). That is, $\mathcal{U}_\x \bigcap \mathcal{H}_{\x, \z_1} = \emptyset$, which implies $\x \in \mathcal{B}_{\z_1}$, and hence $\mathcal{B}_\z \subset \mathcal{B}_{\z_1}$. \vskip 2mm
\medskip

Next, since $\z_1$ lies in the interior of $\mathcal{M}(\mathcal{X}^n)$, we can find a small enough $\varepsilon > 0$ such that $\{\x: \|\x - \z_1\| < \varepsilon\} \subset \mathcal{M}(\mathcal{X}^n)$. For $\forall \u_{\z_1} \in \mathcal{U}_{\z_1}$, if the hyperplane $\{\x: \u_{\z_1}^\top (\x - \z_1) = 0\}$ contains a sample point, say $X_1$, then let $\widetilde{\z}_1 = -\frac{1}{2} \varepsilon \u_{\z_1} + \z_1$. Clearly, $\widetilde{\z}_1 \in \mathcal{M}(\mathcal{X}^n)$, and $\u_{\z_1}^\top \widetilde{\z}_1 < \u_{\z_1}^\top \z_1 = \u_{\z_1}^\top X_1$. As a result, $D(\widetilde{\z}_1, \mathcal{X}^n) \leq P_n(\u_{\z_1}^\top X \leq \u_{\z_1}^\top \widetilde{\z}_1) = P_n(\u_{\z_1}^\top X \leq \u_{\z_1}^\top \z_1) - 1/n = \lambda^*(\mathcal{X}^n) - 1 / n$, contradicting with $\widetilde{\z}_1 \in \mathcal{M}(\mathcal{X}^n)$. \vskip 2mm
\medskip

Hence, $\forall \u_{\z_1} \in \mathcal{U}_{\z_1}$, $\{\x: \u_{\z_1}^\top (\x - \z_1) = 0\}$ contains \emph{no} sample point. This fact implies that there exists a permutation $(i_1, i_2, \cdots, i_n)$ of $(1, 2, \cdots, n)$ satisfying
\begin{eqnarray}
\label{PermU}
    &\u_{\z_1}^\top X_{i_1} \leq \u_{\z_1}^\top X_{i_2} \leq \cdots \leq \u_{\z_1}^\top X_{i_{k^*}} < \u_{\z_1}^\top \z_1 < \nonumber \\
    &\u_{\z_1}^\top X_{i_{k^* + 1}} \leq \cdots \leq \u_{\z_1}^\top X_{i_n},
\end{eqnarray}
where $k^* = n \lambda^*(\mathcal{X}^n)$. Then similar to \cite{LZW2013}, we have that a direction vector $\u$ should belong to $\mathcal{U}_{\z_1}$ if it satisfies
\begin{eqnarray*}
    \left\{
    \begin{array}{ccc}
    \u^\top (X_{i_1} - \z_1) & < & 0\\
    \u^\top (X_{i_2} - \z_1) & < & 0\\
    \vdots &   &  \\
    \u^\top (X_{i_{k^*}} - \z_1) & < & 0\\
    \u^\top (\z_1 - X_{i_{k^*+1}}) & < & 0\\
    \vdots &   &  \\
    \u^\top (\z_1 - X_{i_n}) & < & 0.
    \end{array}
    \right.
\end{eqnarray*}
Using this, it is trivial to check that $\mathcal{U}_{\z_1}$ is non-coplanar when $\z_1$ lies in the interior of $\mathcal{M}(\mathcal{X}^n)$. Hence, there must exist $\vv_{11}, \vv_{12}, \cdots, \vv_{1d} \in \mathcal{U}_{\z_1}$, which are of affine dimension $d$.
\medskip

Observe that, $\forall \z \in \widetilde{\mathcal{B}}_{\z_1}$, there $\exists \widetilde{\vv} \in \{\vv_{11}, \vv_{12}, \cdots, \vv_{1d}\}$ satisfying
\begin{eqnarray}
\label{INCSeq}
  \widetilde{\vv}^\top (\z - \z_1) > 0.
\end{eqnarray}
(If not, $\vv_{1l}^\top (\z - \z_1) = 0$ for $l = 1, 2, \cdots, d$, which lead to $\z = \z_1$. This is impossible due to $\z \in \widetilde{\mathcal{B}}_{\z_1}$.) Hence, $\mathcal{U}_{\z_1} \cap \mathcal{H}_{\z_1, z} \neq \emptyset$, and then $\z_1 \notin \mathcal{B}_\z$. This completes the proof of (\textbf{o3}). \hfill $\Box$
\medskip

Relying on Lemma 2, we are able to find a point $\x_0$ in the interior of $\mathcal{M}(\mathcal{X}^n)$, which lies outside of at least one optimal halfspace of any $\x  \neq \x_0$. Here by optimal halfspace of $x$ we mean the halfspace realizing the depth at $x$. That is, we have the following lemma.
\medskip

\textbf{Lemma 3}. When $\mathcal{X}^n$ is \textbf{IGP}, there must exist an $\x_0 \in \mathcal{M}(\mathcal{X}^n)$ such that $\mathcal{U}_\x \bigcap \mathcal{H}_{\x, \x_0} \neq \emptyset$ for $\forall \x \in \mathcal{R}^d \setminus \{\x_0\}$, i.e., we can find a $\u \in \mathcal{U}_\x$ satisfying $\u^\top \x < \u^\top \x_0$.
\medskip

In the sequel the major task is to prove: $\mathcal{U}_\x \bigcap \mathcal{H}_{\x, \x_0} \neq \emptyset$ for $\forall \x \in \mathcal{M}(\mathcal{X}^n) \setminus \{\x_0\}$ when $d \ge 2$. It consists of three parts, i.e., (\textbf{A}), (\textbf{B}) and (\textbf{C}), related respectively to three scenarios of the affine dimension $\textbf{dim}(\mathcal{M})$ of $\mathcal{M}(\mathcal{X}^n)$. Both (\textbf{A}) and (\textbf{B}) indicates that taking $\x_0 = T^*(\mathcal{X}^n)$ is valid, while (\textbf{C}) is technically much more difficult and the resulted $\x_0$ may $\neq T^*(\mathcal{X}^n)$. In (\textbf{C}), we first obtain a candidate point, say $\bar{\z}_0$, through an iterative procedure consisting of three steps, i.e., (\textbf{a}), (\textbf{b}) and (\textbf{c}), and then show that $\bar{\z}_0$ can serve as $\x_0$. For convenience, we use the same notations (e.g., $\mathcal{X}^n$) as before in Lemma 2 though. Its result can be applied to any other \textbf{IGP} data set, nevertheless.
\medskip

\textbf{Proof}. When $d = 1$, by letting $\x_0$ be the sample median, the proof is trivial. When $d \ge 2$, for $\forall \x \notin \mathcal{M}(\mathcal{X}^n)$ and $\forall \z \in \mathcal{M}(\mathcal{X}^n)$, we claim
$\mathcal{U}_\x \cap \mathcal{H}_{\x, \z} \neq \emptyset$. \emph{If not}, for $\forall \u \in \mathcal{U}_\x$, we have $\u^\top \x \ge \u^\top \z$, which leads to $\lambda^*(\mathcal{X}^n) \leq P_n(\u^\top X \leq \u^\top \z) \leq P_n(\u^\top X \leq \u^\top \x) = D(\x, \mathcal{X}^n)$, contradicting the definition of $\lambda^*$ and $\mathcal{M}(\mathcal{X}^n)$. In the sequel we show that there $\exists \x_0 \in \mathcal{M}(\mathcal{X}^n)$ satisfying $\mathcal{U}_\x \bigcap \mathcal{H}_{\x, \x_0} \neq \emptyset$ for $\forall \x \in \mathcal{M}(\mathcal{X}^n) \setminus \{\x_0\}$.
\medskip

(\textbf{A}) \emph{Scenario $\textbf{dim}(\mathcal{M}) = 0$.} Since $\mathcal{M}(\mathcal{X}^n) = \{T^*(\mathcal{X}^n)\}$, Lemma 3 already holds by letting $\x_0 = T^*(\mathcal{X}^n)$.
\medskip

(\textbf{B}) \emph{Scenario $0 < \textbf{dim}(\mathcal{M}) < d$.} We now show that taking $\x_0 = T^*(\mathcal{X}^n)$ is valid.
\medskip

Relying on Theorem 4.2 of \cite{PS2011}, it is easy to check that there $\exists \mu_* \in \{\mu_j\}_{j=1}^{N_n^d}$ normal to the hyperplane $\Pi_0 = \{\z \in \mathcal{R}^d : \mu_*^\top \z = q_*\}$ such that: (i) $\Pi_0 \supset \mathcal{M}(\mathcal{X}^n)$, (ii) $\Pi_0 $ contains $d$ observations, say $\mathbf{Z}_d := \{X_{k_1}, X_{k_2}, \cdots, X_{k_d}\}$. Here $q_* = \inf\{t \in \mathcal{R}^1: P_n(\mu_*^\top X \leq t) \ge \lambda^*(\mathcal{X}^n)\}$.
\medskip

Obviously, $\mathcal{M}(\mathcal{X}^n) \subset \textbf{cov}(\mathbf{Z}_d)$, i.e., the convex hull of $\mathbf{Z}_d$. \emph{If not}, one may deviate $\Pi_0$ around a point $\x \in \mathcal{M}(\mathcal{X}^n) \setminus \textbf{cov}(\mathbf{Z}_d)$, similar to Theorem 1 of \cite{LLZ2015}, to get rid of $\mathbf{Z}_d$ to obtain a contradiction.
\medskip

For $i = 1, 2, \cdots, d$, let $\mathbf{W}_i$ be the $(d-2)$-dimensional hyperplane passing through $\mathbf{Z}_d \setminus \{X_{k_i}\}$ ($\mathbf{W}_i$ is a singleton when $d=2$), and let $\nu_i \in \mathcal{S}^{d-1}$ be the vector orthogonal to both $\mu_*$ and $\mathbf{W}_i$, and satisfying $\nu_i^\top (X_{k_i} - X_{k_l}) > 0$ for $\forall l \in \{1, 2, \cdots, d\} \setminus \{i\}$. In the following, we show that
\begin{eqnarray}
\label{MXnCovX}
    \mathcal{M}(\mathcal{X}^n) \setminus \{\x_0\} = \bigcup_{i=1}^d \mathcal{D}_i,
\end{eqnarray}
where $\mathcal{D}_i = \left\{\x\in \mathcal{M}(\mathcal{X}^n): \nu_i^\top \x_0 < \nu_i^\top \x \right\}$.
\medskip

\emph{The `$\supset$' part is trivial. We only show the `$\subset$' part.} In fact, if $\exists \z\in \mathcal{M}(\mathcal{X}^n)\setminus \{\x_0\}$ but $\z \notin \bigcup_{i=1}^d \mathcal{D}_i$, then we have
$\nu_i^\top (\x_0 - \z) \ge 0$ for $i = 1, 2, \cdots, d$. Using this and the fact $\mathcal{M}(\mathcal{X}^n) \subset \textbf{cov}(\mathbf{Z}_d)$, we obtain, for $\forall \delta > 0$,
\begin{eqnarray*}
    \delta(\x_0 - \z) + \z &\in& \bigcap_{i=1}^d \left\{\x\in \Pi_0: \nu_i^\top \x \ge \nu_i^\top \z \right\} \\
    &\subset& \bigcap_{i=1}^d \left\{\x\in \Pi_0: \nu_i^\top \x \ge \nu_i^\top X_{k_l}, ~ l \in \{1, 2, \cdots, d\}\setminus \{i\} \right\} =  \textbf{cov}(\mathbf{Z}_d),
\end{eqnarray*}
contradicting with the boundedness of $\textbf{cov}(\mathbf{Z}_d)$. Hence, \eqref{MXnCovX} is true.
\medskip

Relying on \eqref{MXnCovX}, it is easy to find a $\nu \in \{\nu_j\}_{j=1}^d$ and $\varepsilon > 0$ such that $\bar{\mu}_* = \mu_* - \varepsilon \nu$ with $\bar{\mu}_*$ satisfying $P_n(\bar{\mu}_*^\top X \leq \bar{\mu}_*^\top \x) = \lambda^*(\mathcal{X}^n)$ and $\bar{\mu}_*^\top \x < \bar{\mu}_*^\top \x_0$ for $\forall \x \in \mathcal{M}(\mathcal{X}^n) \setminus \{\x_0\}$.
\medskip

(\textbf{C}) \emph{Scenario $\textbf{dim}(\mathcal{M}) = d$.} Let $\z_1 = T^*(\mathcal{X}^n)$. Clearly, $\z_1$ lies in the interior of $\mathcal{M}(\mathcal{X}^n)$. \emph{In the sequel we first obtain a candidate point, say $\bar{\textbf{z}}_0$, through an iterative procedure, and then prove that it can serve as $\bar{\textbf{x}}_0$.}
\medskip

\textbf{Step} (\textbf{a}) If $\mathcal{B}_{\z_1} = \{\z_1\}$, let $\x_0 = \z_1$. This lemma already holds. \emph{Otherwise}, $\widetilde{\mathcal{B}}_{\z_1} \neq \emptyset$, and let
\begin{eqnarray*}
    h(\z_1) = \sup_{\substack{\vv\in \mathcal{U}_{\z_1}, ~ \z \in \widetilde{\mathcal{B}}_{\z_1}}} \vv^\top (\z - \z_1).
\end{eqnarray*}
By the property of the supremum, for $\varepsilon_2 = 1/2$, there must $\exists \z_2 \in \widetilde{\mathcal{B}}_{\z_1}$ and $\exists \bar{\vv}_1 \in \mathcal{U}_{\z_1}$ satisfying
\begin{eqnarray*}
    h(\z_1) - \varepsilon_2 < \bar{\vv}_1^\top (\z_2 - \z_1) \leq h(\z_1).
\end{eqnarray*}

\textbf{Step} (\textbf{b}) Similar to (\textbf{a}), if $\mathcal{B}_{\z_2} = \{\z_2\}$, let $\x_0 = \z_2$ and end the proof of this lemma. \emph{Otherwise}, for $\varepsilon_3 = 1/3$, we similarly have a $\z_3 \in \widetilde{\mathcal{B}}_{\z_2} \subset \mathcal{B}_{\z_1}$ and a $\bar{\vv}_2 \in \mathcal{U}_{\z_2} \subset \mathcal{U}_{\z_1}$, by Lemma 2, satisfying
\begin{eqnarray*}
    h(\z_2) - \varepsilon_3 < \bar{\vv}_2^\top (\z_3 - \z_2) \leq h(\z_2).
\end{eqnarray*}

\textbf{Step} (\textbf{c}) If there is a finite $m$ ($m \ge 3$) such that $\mathcal{B}_{\z_m} = \{\z_m\}$, then let $\x_0 = \z_m$ and end the proof of this lemma. \emph{Otherwise}, by repeating (\textbf{a}) and (\textbf{b}), we can obtain a series of different points $\{\z_i\}_{i=1}^{\infty} \subset \mathcal{M}(\mathcal{X}^n)$ satisfying: for $\forall m > 1$, (\textbf{p1}) $\z_m \in \widetilde{\mathcal{B}}_{\z_{m - 1}}$, (\textbf{p2}) $\mathcal{U}_{\z_{m}} \subset \mathcal{U}_{\z_{m-1}}$, (\textbf{p3}) $\u^\top \z_m \ge \cdots \ge \u^\top \z_2 \ge \u^\top \z_1$ for $\forall \u \in \mathcal{U}_{\z_m}$, (\textbf{p4}) there exists a $\bar{\vv}_{m-1} \in \mathcal{U}_{\z_{m - 1}}$ such that $h(\z_{m-1}) - \varepsilon_m < \bar{\vv}_{m-1}^\top (\z_m - \z_{m-1}) \leq h(\z_{m-1})$, where $\varepsilon_m = 1/m$.
\vskip 2mm

Since $\mathcal{M}(\mathcal{X}^n)$ is bounded, $\{\z_i\}_{i=1}^{\infty}$ must contain a convergent subsequence, say $\{\z_{k_m}\}_{m = 1}^\infty$. Without confusion, suppose $\lim\limits_{m \rightarrow \infty} \z_{k_m} = \bar{\z}_0$. Obviously, $\bar{\z}_0$ should lie \emph{in the interior} of $\mathcal{M}(\mathcal{X}^n)$, because for any point $\x$ on the boundary of $\mathcal{M}(\mathcal{X}^n)$, it is easy to find a $\u_\x \in \mathcal{U}_{\x} \cap \mathcal{H}_{\x, \z}$ for some inner points $\z$ of $\mathcal{M}(\mathcal{X}^n)$.
\vskip 2mm

\textbf{Now we show that $\bar{\z}_0$ \emph{can serve as} $\x_0$.} By (\textbf{p3}), the fact $k_m-1 \ge k_{m-1}$ implies $\u^\top \z_{k_m-1} \ge \u^\top \z_{k_{m-1}}$ for $\forall \u \in \mathcal{U}_{\z_{k_m-1}}$. Hence, for $\bar{\vv}_{m-1} \in \mathcal{U}_{\z_{k_m-1}}$ given in (\textbf{p4}), we have
\begin{eqnarray*}
  h(\z_{k_m-1}) - \varepsilon_{k_m} < \bar{\vv}_{{k_m-1}}^\top (\z_{k_m} - \z_{k_m-1}) \leq \bar{\vv}_{{k_m-1}}^\top (\z_{k_m} - \z_{k_{m-1}}) \leq \|\z_{k_m} - \z_{k_{m-1}}\|.
\end{eqnarray*}
This, together with the convergence of $\{\z_{k_m}\}_{m = 1}^\infty$, leads to
\begin{eqnarray}
\label{hzkm}
  h(\z_{k_m-1}) \rightarrow 0, \quad \text{ as }k_m \rightarrow +\infty.
\end{eqnarray}
\emph{Based on this, we can show that} $\mathcal{B}_{\bar{\z}_0} = \{\bar{\z}_0\}$ through two steps as follows.
\medskip

\textbf{Firstly, for $\forall \bar{k} \in \{k_m\}_{m=1}^\infty$, we have $\bar{\z}_0 \in \mathcal{B}_{\z_{\bar{k}-1}}$.} \emph{If not}, there must exist a $\bar{\u} \in \mathcal{U}_{\bar{\z}_0}$ satisfying $\bar{\u}^\top \bar{\z}_0 < \bar{\u}^\top \z_{\bar{k}-1}$. For $\bar{\u}$, similar to \eqref{PermU}, there exists a permutation $(i_1', i_2', \cdots, i_n')$ of $(1, 2, \cdots, n)$ such that
\begin{eqnarray*}
    \bar{\u}^\top X_{i_1'} \leq \bar{\u}^\top X_{i_2'} \leq \cdots \leq \bar{\u}^\top X_{i_{k^*}'} < \bar{\u}^\top \bar{\z}_0 < \bar{\u}^\top X_{i_{k^*+1}'} \leq \cdots \leq \bar{\u}^\top X_{i_n'}.
\end{eqnarray*}
Let $\delta_0 = \frac{1}{2} \min \left\{\bar{\u}^\top (\bar{\z}_0 - X_{i_{k^*}'}),~ \bar{\u}^\top ( X_{i_{k^*+1}'} - \bar{\z}_0), ~ \bar{\u}^\top (\z_{\bar{k}-1} - \bar{\z}_0)\right\}$. By the convergence of $\{\z_{k_m}\}_{m = 1}^\infty$, we can find a $k_m' > \bar{k}$ among $\{k_m\}_{m=1}^\infty$ such that $\| \z_{k_m'} - \bar{\z}_0 \| < \delta_0$. This, combined with $\left|\bar{\u}^\top (\z_{k_m'} - \bar{\z}_0) \right| \leq \| \z_{k_m'} - \bar{\z}_0 \|$, leads to
\begin{eqnarray}
\label{UKM}
    \bar{\u}^\top X_{i_1'} \leq \bar{\u}^\top X_{i_2'} \leq \cdots \leq \bar{\u}^\top X_{i_{k^*}'} < \bar{\u}^\top \z_{k_m'} < \bar{\u}^\top X_{i_{k^*+1}'} \leq \cdots \leq \bar{\u}^\top X_{i_n'}.
\end{eqnarray}
That is, $\bar{\u} \in \mathcal{U}_{\z_{k_m'}}$ ($\subset \mathcal{U}_{\z_{\bar{k}-1}}$). On the other hand, following a similar fashion to \eqref{UKM}, we have
\begin{eqnarray*}
    \bar{\u}^\top \z_{k_m'} < \bar{\u}^\top \z_{\bar{k}-1}.
\end{eqnarray*}
This nevertheless contradicts with (\textbf{o1}) and (\textbf{o2}) of Lemma 2 due to $\z_{k_m'} \in \mathcal{B}_{\z_{\bar{k}-1}}$ when $k_m' > \bar{k}$. Hence, $\bar{\z}_0 \in \mathcal{B}_{\z_{\bar{k}-1}}$.
\medskip

\textbf{Secondly, we show $\mathcal{B}_{\bar{\z}_0} \setminus \{\bar{\z}_0\} = \emptyset$ based on the fact $\bar{\z}_0 \in \mathcal{B}_{\z_{\bar{k}-1}}$ for $\forall \bar{k} \in \{k_m\}_{m=1}^\infty$.} \emph{If not}, suppose $\x \in \mathcal{B}_{\bar{\z}_0} \setminus \{\bar{\z}_0\}$ without loss of generality. Then, similar to \eqref{INCSeq}, we can find a $\vv_0 \in \mathcal{U}_\x$ such that $\vv_0^\top \x - \vv_0^\top \bar{\z}_0 > 0$. On the other hand, by (\textbf{o2})-(\textbf{o3}) of Lemma 2, the facts $\x \in \mathcal{B}_{\bar{\z}_0}$ and $\bar{\z}_0 \in \mathcal{B}_{\z_{\bar{k}-1}}$ together imply $\x \in \mathcal{B}_{\z_{\bar{k}-1}} \setminus \{\z_{\bar{k}-1}\}$ and $\vv_0 \in \mathcal{U}_{\z_{\bar{k}-1}}$. These, combined with (\textbf{o1}) of Lemma 2 and the property of the supremum, lead to
\begin{eqnarray*}
  h(\z_{\bar{k}-1}) \ge \vv_0^\top \x - \vv_0^\top \z_{\bar{k}-1} \ge \vv_0^\top \x - \vv_0^\top \bar{\z}_0 > 0.
\end{eqnarray*}
Nevertheless, $\vv_0^\top \x - \vv_0^\top \bar{\z}_0$ does \emph{not} depend on $\bar{k}$, contradicting with \eqref{hzkm}.
\medskip

This completes the proof of this lemma. \hfill $\Box$

\vskip 0.1 in
\section{FSBP of Tukey's halfspace median (Main results)}\label{main}
\paragraph{}
\vskip 0.1 in \label{UpperBound}

Note that for $\u \in \mathcal{S}^{d-1}$, its $\mathbb{A}_\u$-projections $\mathbf{X}_\u^n$ is \emph{not} \textbf{IGP} if \eqref{eqnu0} is violated, while in the proof of our main theorem, we have to handle such situations that $\mathbf{X}_\u^n$ is not \textbf{IGP}. Hence, in addition to three preliminary lemmas above, we need three more lemmas as follows.
\medskip

\textbf{Lemma 4}. There exists a $\u_0 \in \mathcal{S}^{d-1}$ such that: (\textbf{s1}) $\u_0$ satisfies \eqref{eqnu0}, and (\textbf{s2})
\begin{eqnarray}
\label{nv0}
  \lambda^*(\mathbf{X}_{\u_0}^n) = \inf_{\u \in \mathcal{S}^{d-1}} \lambda^*(\mathbf{X}_\u^n).
\end{eqnarray}

\textbf{Proof}. Note that $\lambda^*(\mathbf{X}_\u^n) \in \{0, 1/n, \cdots, 1\}$ for $\forall \u \in \mathcal{S}^{d-1}$. Hence, there must exist a $\u_0 \in \mathcal{S}^{d-1}$ satisfying \eqref{nv0}. \emph{This completes the proof of} (\textbf{s2}).
\medskip

\emph{Now we show} (\textbf{s1}). For simplicity, let $\mathcal{N}_0 = \{1, 2, 3, \cdots \}$ and denote $\inf_{k \in \mathcal{I}} t_k$ as the infimum of the set $\{t: t = t_k, k \in \mathcal{I}\}$ related to $\{t_k\}_{k = 1}^\infty \subset \mathcal{R}^1$, where $\mathcal{I}$ denotes some subscript sets that $\mathcal{I} \subset \mathcal{N}_0$.
\medskip

By noting that $\mathcal{S}^{d-1}$ is of affine dimension $d$, it is easy to check that $N_n^d$ hyperplanes $\Pi_j = \{\x \in \mathcal{R}^d: \mu_j^\top \x = 0\}$, $j = 1, 2, \cdots, N_n^d$, together divide $\mathcal{S}^{d-1}$ into only a finite number of non-coplanar fragments. Hence, for any $\vartheta_0 \in \mathcal{S}^{d-1}$ such that $\mu_j^\top \vartheta_0 = 0$ for some $j \in \{1, 2, \cdots, N_n^d\}$, we can find a sequence $\{\vartheta_k\}_{k = 1}^\infty \subset \mathcal{S}^{d-1}$ satisfying: (i) each $\vartheta_k$ lies in the interior of a non-coplanar fragment of $\mathcal{S}^{d-1}$ and satisfies display \eqref{eqnu0}, and (ii) $\lim\limits_{k \rightarrow \infty} \vartheta_k = \vartheta_0$. Now we show that
\begin{eqnarray*}
  \lambda^*(\mathbf{X}_{\vartheta_0}^n) \ge \inf_{k \in \mathcal{N}_0} \lambda^*(\mathbf{X}_{\vartheta_k}^n).
\end{eqnarray*}

Since $\vartheta_k$ can be obtained through rotating $\vartheta_0$, there must exist a unique orthogonal matrix $\mathbb{Q}_k$ such that $\vartheta_k = \mathbb{Q}_k \vartheta_0$. Obviously, $\lim\limits_{k \rightarrow \infty} \vartheta_k = \vartheta_0$ implies $\lim\limits_{k \rightarrow \infty} \mathbb{Q}_k = \mathbb{I}_p$, which is the $d\times d$ identical matrix. Denote $\mathbb{A}_0 := \mathbb{A}_{\vartheta_0}$. For $\vartheta_k$, since $\mathbb{Q}_k \mathbb{A}_0$ satisfies display \eqref{mathbbA}, it can serve as $\mathbb{A}_{\vartheta_k}$. For simplicity, hereafter denote $\mathbb{A}_k := \mathbb{A}_{\vartheta_k} = \mathbb{Q}_k \mathbb{A}_0$, and $\theta_k = T^*(\mathbf{X}_{\vartheta_k}^n)$ for $k \in \mathcal{N}_0$. Note that, for any $u \in \mathcal{S}^{d-1}$, $|\u^\top X_i| \leq \max\limits_{1\leq j \leq n} \|X_j\|$. Hence, $\{\theta_k\}_{k=1}^\infty$ is bounded, and it should contain a convergent subsequence. Without confusion, suppose $\{\theta_k\}_{k=1}^\infty$ is convergent with $\lim\limits_{k \rightarrow \infty} \theta_k = \theta_0$. (If not, use the convergent subsequence as $\{\theta_k\}_{k=1}^\infty$ instead).
\medskip

Suppose $\textbf{v}_0 \in \mathcal{S}^{d-2}$ satisfies that $\mathbf{p}_n (\textbf{v}_0^\top \mathbf{X}_0 \leq \textbf{v}_0^\top \theta_0) = D(\theta_0, \mathbf{X}_{\vartheta_0}^n)$, where $\mathbf{X}_0 = \mathbb{A}_0^\top X$, and hereafter $\mathbf{p}_n$ denotes the empirical probability measure in the $(d-1)$-dimensional space. For convenience, let
\begin{eqnarray*}
  \mathcal{J}^0 &=& \{j : \textbf{v}_0^\top (\mathbb{A}_0^\top X_j) = \textbf{v}_0^\top \theta_0, ~ j = 1, 2, \cdots, n\},\\
  \mathcal{J}^- &=& \{j : \textbf{v}_0^\top (\mathbb{A}_0^\top X_j) < \textbf{v}_0^\top \theta_0, ~ j = 1, 2, \cdots, n\},\\
  \mathcal{J}^+ &=& \{j : \textbf{v}_0^\top (\mathbb{A}_0^\top X_j) > \textbf{v}_0^\top \theta_0, ~ j = 1, 2, \cdots, n\}.
\end{eqnarray*}
Obviously, (i) $n \mathbf{p}_n (\textbf{v}_0^\top \mathbf{X}_0 \leq \textbf{v}_0^\top \theta_0) = \#(\mathcal{J}^0 \cup \mathcal{J}^-)$, where $\# (\mathcal{A})$ denotes the cardinal number of a set $\mathcal{A}$, and (ii) $\# (\mathcal{J}^0) \leq d$ when $\mathcal{X}^n$ is \textbf{IGP}. Without loss of generality, write $\mathcal{J}^0 = \{j_1, j_2, \cdots, j_q\}$, where $0 \leq q \leq d$.
\begin{enumerate}
  \item[(i)] When $q = 0$, i.e., $\mathcal{J}^0 = \emptyset$, the facts $\lim\limits_{k \rightarrow \infty} \mathbb{Q}_k = \mathbb{I}_p$ and $\lim\limits_{k \rightarrow \infty} \theta_k = \theta_0$ together lead to
      \begin{eqnarray}
      \label{eqnA0Theta0}
          \lim\limits_{k \rightarrow \infty} I\left(\textbf{v}_0^\top (\mathbb{A}_k^\top X_i) \leq \textbf{v}_0^\top \theta_k \right) &=& \lim\limits_{k \rightarrow \infty} I\left(\textbf{v}_0^\top (\mathbb{A}_0^\top \mathbb{Q}_k^\top X_i) \leq \textbf{v}_0^\top \theta_k \right)\nonumber\\
          &=& I\left(\textbf{v}_0^\top (\mathbb{A}_0^\top X_i) \leq \textbf{v}_0^\top \theta_0 \right)
      \end{eqnarray}
      for each $i = 1, 2, \cdots, n$, where $I(\cdot)$ denotes the indicative function. Using this, we obtain
      \begin{eqnarray}
      \label{PN1}
        \lim\limits_{k \rightarrow \infty} \mathbf{p}_n (\textbf{v}_0^\top \mathbf{X}_k \leq \textbf{v}_0^\top \theta_k) = \mathbf{p}_n (\textbf{v}_0^\top \mathbf{X}_0 \leq \textbf{v}_0^\top \theta_0),
      \end{eqnarray}
      where $\mathbf{X}_k = \mathbb{A}_k^\top X$.
      \medskip

  \item[(ii)] When $q > 0$, i.e., $\mathcal{J}^0 \neq \emptyset$, we have the following results. \medskip

      For $l = 1, 2, \cdots, q$, denote
      \begin{eqnarray*}
        \mathcal{N}_{l - 1}^- &=& \left\{k \in \mathcal{N}_{l - 1}: \textbf{v}_0^\top (\mathbb{A}_k^\top X_{j_l}) \leq \textbf{v}_0^\top \theta_k \right\},\\
        \mathcal{N}_{l - 1}^+ &=& \left\{k \in \mathcal{N}_{l - 1}: \textbf{v}_0^\top (\mathbb{A}_k^\top X_{j_l}) > \textbf{v}_0^\top \theta_k \right\}.
      \end{eqnarray*}

      Check whether or not $\# (\mathcal{N}_{l - 1}^-) < \infty$. If not, set $\mathcal{N}_l = \mathcal{N}_{l-1}^-$; Otherwise, $\mathcal{N}_l = \mathcal{N}_{l - 1}^+$. Then $\max\{\#(\mathcal{N}_{l - 1}^-), ~\#(\mathcal{N}_{l - 1}^-)\} = \infty$ due to $\#(\mathcal{N}_{l-1}) = \infty$ for $l = 1, 2, \cdots, q$. Hence, $\#(\mathcal{N}_q) = \infty$ because $q \leq d$.
      \medskip

      Using this, we claim that for each $l = 1, 2, \cdots, q$, either
      \begin{eqnarray}
      \label{case1}
        I\left(\textbf{v}_0^\top (\mathbb{A}_k^\top X_{j_l}) \leq \textbf{v}_0^\top \theta_k\right) = 0, \quad \text{for all } k \in \mathcal{N}_q,
      \end{eqnarray}
      or
      \begin{eqnarray}
      \label{case2}
        I\left(\textbf{v}_0^\top (\mathbb{A}_k^\top X_{j_l}) \leq \textbf{v}_0^\top \theta_k\right) = 1, \quad \text{for all } k \in \mathcal{N}_q
      \end{eqnarray}
      is true by the construction of $\mathcal{N}_q$. Hence,
      \begin{eqnarray*}
        \lim_{k \in \mathcal{N}_q,~ k \rightarrow \infty} I\left(\textbf{v}_0^\top (\mathbb{A}_k^\top X_{j_l}) \leq \textbf{v}_0^\top \theta_k\right) = \left\{
        \begin{array}{l}
            0,  \quad \text{if \eqref{case1} is true}\\ \\
            1,  \quad \text{if \eqref{case2} is true}.
        \end{array}
         \right.
      \end{eqnarray*}
      This, together with $\textbf{v}_0^\top (\mathbb{A}_0^\top X_i) = \textbf{v}_0^\top \theta_0$ ($i \in \mathcal{J}^0$), leads to
      \begin{eqnarray}
      \label{J0}
        \lim_{k \in \mathcal{N}_q,~ k \rightarrow \infty} I\left(\textbf{v}_0^\top (\mathbb{A}_k^\top X_i) \leq \textbf{v}_0^\top \theta_k\right) \leq I\left(\textbf{v}_0^\top (\mathbb{A}_0^\top X_i) \leq \textbf{v}_0^\top \theta_0\right) = 1,
      \end{eqnarray}
      for $\forall i \in \mathcal{J}^0$. \medskip

      On the other hand, similar to \eqref{eqnA0Theta0}, we have
      \begin{eqnarray*}
          \lim_{k \in \mathcal{N}_q,~ k \rightarrow \infty} I\left(\textbf{v}_0^\top (\mathbb{A}_k^\top X_i) \leq \textbf{v}_0^\top \theta_k \right) = I\left(\textbf{v}_0^\top (\mathbb{A}_0^\top X_i) \leq \textbf{v}_0^\top \theta_0 \right)
      \end{eqnarray*}
       for $\forall i \notin \mathcal{J}^0$. This, together with \eqref{J0}, shows
       \begin{eqnarray}
       \label{PN2}
        \lim_{k \in \mathcal{N}_q,~ k \rightarrow \infty} \mathbf{p}_n (\textbf{v}_0^\top \mathbf{X}_k \leq \textbf{v}_0^\top \theta_k) \leq \mathbf{p}_n (\textbf{v}_0^\top \mathbf{X}_0 \leq \textbf{v}_0^\top \theta_0).
      \end{eqnarray}
\end{enumerate}

Next, by observing
\begin{eqnarray*}
  \mathbf{p}_n (\textbf{v}_0^\top \mathbf{X}_k \leq \textbf{v}_0^\top \theta_k) \ge \inf_{\textbf{v} \in\mathcal{S}^{d-2}} \mathbf{p}_n (\textbf{v}^\top \mathbf{X}_k \leq \textbf{v}^\top \theta_k) = D(\theta_k, \mathbf{X}_{\vartheta_k}^n) = \lambda^*(\mathbf{X}_{\vartheta_k}^n), ~k = 1, 2, \cdots,
\end{eqnarray*}
and the fact that $\mathbf{p}_n(\cdot) \in \{0, 1/n , 2/n, \cdots, 1\}$, we have that
\begin{eqnarray*}
  \lim\limits_{k \in \widetilde{\mathcal{I}}, ~ k \rightarrow \infty} \mathbf{p}_n (\textbf{v}_0^\top \mathbf{X}_k \leq \textbf{v}_0^\top \theta_k) \ge \inf_{l \in \widetilde{\mathcal{I}}} \lambda^*(\mathbf{X}_{\vartheta_l}^n) \ge \inf_{j \in \mathcal{N}_0} \lambda^*(\mathbf{X}_{\vartheta_j}^n),
\end{eqnarray*}
where $\widetilde{\mathcal{I}} = \mathcal{N}_0$ if $\mathcal{J}^0 = \emptyset$, otherwise $\widetilde{\mathcal{I}} = \mathcal{N}_q$. This, combined with \eqref{PN1} and \eqref{PN2}, implies
\begin{eqnarray*}
  \lambda^*(\mathbf{X}_{\vartheta_0}^n) \ge D(\theta_0, \mathbf{X}_{\vartheta_0}^n) = \mathbf{p}_n (\textbf{v}_0^\top \mathbf{X}_0 \leq \textbf{v}_0^\top \theta_0) \ge \inf_{k \in \mathcal{N}_0} \lambda^*(\mathbf{X}_{\vartheta_k}^n).
\end{eqnarray*}
Finally, by noting that the image of $\lambda^*(\mathbf{X}_{\vartheta_k}^n)$ takes only a finite set of values, we claim that there must exist a $k_0 > 0$ such that $\lambda^*(\mathbf{X}_{\vartheta_{k_0}}^n) = \inf_{k \in \mathcal{N}_0} \lambda^*(\mathbf{X}_{\vartheta_k}^n)$. This lemma then follows immediately.  \hfill $\Box$
\medskip

The aforementioned four lemmas are important in proving the upper bound parts of the main theorem, while the following two lemmas play a key role in obtaining the lower bound of the FSBP of \emph{HM}.
\medskip

\textbf{Lemma 5}. For any given $\y \in \mathcal{R}^d \setminus \textbf{cov}(\mathcal{X}^n)$ ($d \ge 2$), there exists a $\u_0 \in \mathcal{S}^{d-1}$ such that
\begin{eqnarray*}
  D(\mathbb{A}_{\u_0}^\top \y, \mathbf{X}_{\u_0}^n) = \lambda^*(\mathbf{X}_{\u_0}^n).
\end{eqnarray*}

\textbf{Proof}. For $\forall \y \in \mathcal{R}^d \setminus \textbf{cov}(\mathcal{X}^n)$, $\|X_i - \y\| \neq 0$ for $i = 1, 2, \cdots, n$. Hence, we may let $\mathcal{W}^n = \{W_1$, $W_2, \cdots, W_n\}$, where $W_i = (X_i - \y) / \|X_i - \y\|$.
\medskip

Next, by observing the facts that (i) $I(\u^\top X_i \leq \u^\top y) = I(\u^\top W_i \leq 0)$, and (ii) $\textbf{v}^\top \mathbb{A}_\u^\top \u = 0$ and $\|\mathbb{A}_\u \textbf{v}\| = 1$ hold true for $\forall \u \in \mathcal{S}^{d-1}$ and $\forall \textbf{v} \in \mathcal{S}^{d-2}$, we obtain
\begin{eqnarray}
\label{ATD}
  D(\mathbb{A}_\u^\top \y, \mathbf{X}_\u^n) &=& \inf_{\textbf{v} \in \mathcal{S}^{d-2}} \mathbf{p}_n\left(\textbf{v}^\top \mathbf{X} \leq \textbf{v}^\top (\mathbb{A}_\u^\top \y)\right)\nonumber\\
  &=& \inf_{\bar{\u} \in \mathcal{S}^{d-1},~ \bar{\u} \bot \u} P_n\left(\bar{\u}^\top X \leq \bar{\u}^\top \y\right)\nonumber\\
  &=& \inf_{\bar{\u} \in \mathcal{S}^{d-1},~ \bar{\u} \bot \u} \frac{1}{n}\sum_{i=1}^n I(\bar{\u}^\top W_i \leq 0),
\end{eqnarray}
where by $\alpha \bot \beta$ we mean that $\alpha$ is normal to $\beta$ hereafter.
\medskip

Note that $\mathcal{W}^n \in \mathcal{S}^{d-1}$ and $\u$ belongs to the closed hemisphere $\{\vv \in \mathcal{S}^{d-1}: \bar{\u}^\top \vv \leq 0\}$. According to \cite{LS1992}, \eqref{ATD} is in fact the \emph{angular Tukey's depth} of $\u$ with respect to $\mathcal{W}^n$ on the sphere $\mathcal{S}^{d-1}$. Let $\u_0$ be the corresponding angular Tukey's median of $\mathcal{W}^n$. Then this lemma follows immediately.  \hfill $\Box$
\medskip

\textbf{Lemma 6}. Let $\textbf{B}(\mathcal{X}^n)$ be the boundary of the convex hull $\textbf{cov}(\mathcal{X}^n)$ of $\mathcal{X}^n$. Then for any given $\u_\ell \in \mathcal{S}^{d-1}$ ($d \ge 2$), we have that
\begin{eqnarray*}
  D(\z, \mathcal{X}^n \cup \mathcal{Y}^m) \ge \frac{\min\{n \lambda^*(\mathbf{X}_{\u_\ell}^n), ~ m + 1\}}{n + m},
\end{eqnarray*}
where $\mathcal{Y}^m$ denotes the data set containing exactly $m$ repetitions of $\y$ with $\y \in \ell \setminus \textbf{cov}(\mathcal{X}^n)$, and $\z$ the closer to $\y$ intersection of $\ell$ and $\textbf{B}(\mathcal{X}^n)$, where $\ell = \{\x: \x = \mathbb{A}_{\u_\ell} \mathbf{x} + \delta \u_\ell, \delta \in \mathcal{R}^1\}$ with $\mathbf{x} \in \mathcal{M}(\mathbf{X}_{\u_\ell}^n)$.
\medskip

\textbf{Proof}. For any $\u \in \mathcal{S}^{d-1}$, we have the following results.
\begin{enumerate}
  \item[]
    \begin{enumerate}
      \item[(i)]  If $\u^\top \z \ge \u^\top \y$, by observing that $\z \in \mathbf{cov}(\mathcal{X}^n)$ implies $nP_n(\u^\top X \leq \u^\top \z) \ge 1$, we have $(n + m)P_{n+m}(\u^\top X \leq \u^\top \z) \ge m + 1$, where $P_{n+m}$ denotes the empirical probability measure related to the data set $\mathcal{X}^n \cup \mathcal{Y}^m$.

      \item[(ii)] If $\u^\top \z < \u^\top \y$, it is trivial that $(n + m)P_{n+m}(\u^\top X \leq \u^\top \z) =  nP_n(\u^\top X \leq \u^\top \z)$. Now we prove that if there exist a $\vv_0 \in \mathcal{S}^{d-1}$ such that $P_n(\vv_0^\top X \leq \vv_0^\top \z) < \lambda^*(\mathbf{X}_{\u_\ell}^n)$, we can obtain a contradiction.
          \medskip

          Without confusion, let $\textbf{n}_\ell \in \mathcal{S}^{d-1}$ be a normal vector of $\ell$. Denote $\mathcal{Q}_{\textbf{n}_\ell}^1 = \{\x: \textbf{n}_\ell^\top (\x - \z) > 0,~ \vv_0^\top (\x - \z) > 0\}$, $\mathcal{Q}_{\textbf{n}_\ell}^2 = \{\x: \textbf{n}_\ell^\top (\x - \z) < 0,~ \vv_0^\top (\x - \z) > 0\}$, $\mathcal{Q}_{\textbf{n}_\ell}^3 = \{\x: \textbf{n}_\ell^\top (\x - \z) < 0,~ \vv_0^\top (\x - \z) < 0\}$ and $\mathcal{Q}_{\textbf{n}_\ell}^4 = \{\x: \textbf{n}_\ell^\top (\x - \z) > 0,~ \vv_0^\top (\x - \z) < 0\}$.
          \medskip

          Clearly, $\textbf{n}_\ell \neq \pm \vv_0$ when $P_n(\vv_0^\top X \leq \vv_0^\top \z) < \lambda^*(\mathbf{X}_{\u_\ell}^n)$ because of
          \begin{eqnarray*}
            \lambda^*(\mathbf{X}_{\u_\ell}^n) = \inf_{\vv\in \mathcal{S}^{d-1}, \vv \bot \ell} P_n(\vv^\top X \leq \vv^\top \z).
          \end{eqnarray*}
          Among all normal vectors of $\ell$, there must exist at least one $\textbf{n}_\ell$ satisfying (\textbf{c1}): $\min\{\aleph(\mathcal{Q}_{\textbf{n}_\ell}^1),\ \aleph(\mathcal{Q}_{\textbf{n}_\ell}^2)\} = 0$ with $\aleph(\mathcal{A}) = \sum_{i=1}^n I(X_i \in \mathcal{A})$ for a set $\mathcal{A}$. If not, there will exist a contradiction with the facts that $\y \notin \textbf{cov}(\mathcal{X}^n)$ and $\z$ is the closer intersection of $\ell$ and $\textbf{B}(\mathcal{X}^n)$ to $\y$.
          \medskip

          Without loss of generality, suppose $\aleph(\mathcal{Q}_{\textbf{n}_\ell}^2) = 0$. Then, (\textbf{c1}), together with the fact $\aleph(\ell \cap \{\x: \vv_0^\top (\x - \z) > 0\}) = 0$, easily leads to $\aleph(\bar{\mathcal{Q}}_{\textbf{n}_\ell}^2 \cup \bar{\mathcal{Q}}_{\textbf{n}_\ell}^3) = \aleph(\bar{\mathcal{Q}}_{\textbf{n}_\ell}^3)$, where $\bar{\mathcal{A}}$ denotes the closure of $\mathcal{A}$. Note that: $\aleph(\bar{\mathcal{Q}}_{\textbf{n}_\ell}^2 \cup \bar{\mathcal{Q}}_{\textbf{n}_\ell}^3) = n P_n(\textbf{n}_\ell^\top X \leq \textbf{n}_\ell^\top \z) \ge n \lambda^*(\mathbf{X}_{\u_\ell}^n)$. Hence, $\aleph(\bar{\mathcal{Q}}_{\textbf{n}_\ell}^3) \ge n \lambda^*(\mathbf{X}_{\u_\ell}^n)$. Obviously, this contradicts with the assumption such that $n \lambda^*(\mathbf{X}_{\u_\ell}^n) > n P_n(\vv_0^\top X \leq \vv_0^\top \z) = \aleph(\bar{\mathcal{Q}}_{\textbf{n}_\ell}^3 \cup \bar{\mathcal{Q}}_{\textbf{n}_\ell}^4) \ge \aleph(\bar{\mathcal{Q}}_{\textbf{n}_\ell}^3)$.
    \end{enumerate}
\end{enumerate}

Combined with (i) and (ii), we obtain this lemma immediately.  \hfill $\Box$
\medskip

With Lemmas 1-6 at hand, we now are able to prove our main theorem as follows, in which we obtain a precise result on the FSBP for \emph{HM}.
\medskip

\textbf{Theorem 1}. Suppose $\mathcal{X}^n$ are in general position. When $d \ge 2$, the FSBP of Tukey's halfspace median $T^*$ is
\begin{eqnarray*}
  \varepsilon(T^*, \mathcal{X}^n) = \frac{\inf_{\u \in \mathcal{S}^{d-1}} \lambda^*(\mathbf{X}_\u^n)}{1 + \inf_{\u \in \mathcal{S}^{d-1}} \lambda^*(\mathbf{X}_\u^n)}.
\end{eqnarray*}

\textbf{Proof}. Let $\y$ be an arbitrary datum, and assume that $\mathcal{Y}^m$ contains exactly $m$ repetitions of $\y$. Clearly, for $\mathcal{X}^n$ and $\mathcal{X}^n \cup \mathcal{Y}^m$, \textbf{cov}($\mathcal{X}^n) \subset \textbf{cov}(\mathcal{X}^n \cup \mathcal{Y}^m$).
\medskip

By Lemma 4, there is a $\u_0 \in \mathcal{S}^{d-1}$ satisfying display \eqref{eqnu0}, and simultaneously
\begin{eqnarray*}
  \lambda^*(\mathbf{X}_{\u_0}^n) = \inf_{\u \in \mathcal{S}^{d-1}} \lambda^*(\mathbf{X}_\u^n).
\end{eqnarray*}
By Lemma 1, $\mathbf{X}_{\u_0}^n$ is still in general position under the current assumptions. This, combined with Lemma 2, indicates that there $\exists \mathbf{x}_0 \in \mathcal{M}(\mathbf{X}_{\u_0}^n)$ such that: there $\exists \mathbf{u} \in \mathcal{U}_{\mathbf{x}}$ satisfying $\mathbf{u}^\top \mathbf{x} < \mathbf{u}^\top \mathbf{x}_0$ for $\forall \mathbf{x} \neq \mathbf{x}_0$.
\medskip

Let $\x_0 = \mathbb{A}_{\u_0} \mathbf{x}_0$, and $\ell_0 = \{\x: \x = \x_0 + \delta \u_0, \delta \in \mathcal{R}^1\}$. Obviously, for any $\x \in \ell_0$ and $\delta \in \mathcal{R}^1$, we have $\mathbb{A}_{\u_0}^\top \x = \mathbb{A}_{\u_0}^\top \mathbb{A}_{\u_0} \mathbf{x}_0 + \delta \mathbb{A}_{\u_0}^\top \u_0 = \mathbf{x}_0$. That is, the $\mathbb{A}_{\u_0}$-projection of any $\x \in \ell_0$ is $\mathbf{x}_0$.
\medskip

As $\y$ is arbitrary, we suppose $\y \in \ell_0 \setminus \textbf{cov}(\mathcal{X}^n)$. \textbf{Now we show that $n\lambda^*(\mathbf{X}_{\u_0}^n)$ such $\y$ suffice for breaking down $T^*$}.
\medskip

Decompose $\textbf{cov}(\mathcal{X}^n) = \mathcal{D}_1 \bigcup \mathcal{D}_2$, where $\mathcal{D}_1 = \ell_0 \bigcap \textbf{cov}(\mathcal{X}^n)$ and $\mathcal{D}_2 = \textbf{cov}(\mathcal{X}^n) \setminus \ell_0$.
\begin{enumerate}
\item[]
\begin{enumerate}
    \item[(\textbf{i})]  For $\forall \x \in \mathcal{D}_1$, its $\mathbb{A}_{\u_0}$-projection is $\mathbf{x}_0 \in \mathcal{R}^{d-1}$. Hence, there $\exists \mathbf{v} \in \mathcal{U}_{\mathbf{x}_0}$ satisfying $\mathbf{p}_n (\mathbf{v}^\top \mathbf{X} \leq \mathbf{v}^\top \mathbf{x}_0) = \lambda^*(\mathbf{X}_{\u_0}^n)$, where $\mathbf{X} = (\mathbb{A}_{\u_0}^\top X)$. \medskip

        Next, for $\mathbf{v}$, similar to \cite{DM2014}, by making $\varepsilon > 0$ small enough, we have that $\bar{\mathbf{v}} = \mathbf{v} - \varepsilon \mathbb{A}_{u_0}^\top (y - x)$ still satisfies $\lambda^*(\mathbf{X}_{\u_0}^n) = \mathbf{p}_n (\bar{\mathbf{v}}^\top \mathbf{X} \leq \bar{\mathbf{v}}^\top \mathbf{x}_0)$. Using this and the fact that
        $\bar{\u}^\top (\y - \x) = \mathbf{v}^\top (\mathbb{A}_{\u_0}^\top (\y - \x)) - \varepsilon \|\mathbb{A}_{\u_0}^\top (\y - \x)\| = - \varepsilon \|\mathbb{A}_{\u_0}^\top (\y - \x)\| < 0$,
        where $\bar{\u} = \mathbb{A}_{\u_0} \bar{\mathbf{v}}$, we obtain
        \begin{eqnarray*}
          D(\x, \mathcal{X}^n \cup \mathcal{Y}^m) &\leq& \frac{n}{n + m} P_n\left(\bar{\u}^\top X \leq \bar{\u}^\top \x\right)\\
          &=& \frac{n}{n + m} P_n\left(\bar{\mathbf{v}}^\top (\mathbb{A}_{\u_0}^\top X) \leq \bar{\mathbf{v}}^\top (\mathbb{A}_{\u_0}^\top \x)\right)\\
          &=& \frac{n}{n + m} \mathbf{p}_n \left(\bar{\mathbf{v}}^\top \mathbf{X} \leq \bar{\mathbf{v}}^\top \mathbf{x}\right) \\
          &=& \frac{n}{n + m} \lambda^*(\mathbf{X}_{\u_0}^n).
        \end{eqnarray*}

    \item[(\textbf{ii})] For $\forall \x \in \mathcal{D}_2$, denote $\mathbf{x} = \mathbb{A}_{\u_0}^\top \x$. Since $\mathbf{x} \neq \mathbf{x}_0$, by Lemma 2, we can find a $\mathbf{v} \in \mathcal{U}_{\mathbf{x}}$ such that $\mathbf{v}^\top \mathbf{x} < \mathbf{v}^\top \mathbf{x}_0$. By noting $\mathbb{A}_{\u_0}^\top \y = \mathbf{x}_0$, we have $\mathbf{v}^\top \mathbf{x} < \mathbf{v}^\top (\mathbb{A}_{\u_0}^\top \y)$. Using this, a similar derivation to \textbf{(i)} leads to
        \begin{eqnarray*}
          D(\x, \mathcal{X}^n \cup \mathcal{Y}^m) \leq \frac{n}{n + m} P_n\left(\u^\top X \leq \u^\top \x\right) \leq \frac{n}{n + m} \lambda^*(\mathbf{X}_{\u_0}^n),
        \end{eqnarray*}
        where $\u = \mathbb{A}_{\u_0} \textbf{v}$.
    \end{enumerate}
\end{enumerate}

(\textbf{i}) and (\textbf{ii}) lead to
\begin{eqnarray}
\label{inConX}
    \sup_{\x \in \textbf{cov}(\mathcal{X}^n)} D(\x, \mathcal{X}^n \cup \mathcal{Y}^m) \leq \frac{n}{n + m} \lambda^*(\mathbf{X}_{\u_0}^n).
\end{eqnarray}

Next, for any $\x \notin \textbf{cov}(\mathcal{X}^n)$, there must exist a $\u_\x$ such that $P_n(\u_\x^\top X \leq \u_\x^\top \x) = 0$ by the convexity of $\textbf{cov}(\mathcal{X}^n)$. Using this, we claim that
\begin{eqnarray*}
    \u_\x^\top \y \leq \u_\x^\top \x < \u_\x^\top X_1, \u_\x^\top X_2, \cdots, \u_\x^\top X_n
\end{eqnarray*}
hold true for any $\x \in \textbf{cov}(\mathcal{X}^n \cup \mathcal{Y}^m) \setminus \textbf{cov}(\mathcal{X}^n)$ but $\x \neq \y$.
(The fact that $\u_\x^\top \x < \u_\x^\top \y$, $\u_\x^\top X_1$, $\u_\x^\top X_2$, $\cdots, \u_\x^\top X_n$ contradicts with the convexity of $\textbf{cov}(\mathcal{X}^n \cup \mathcal{Y}^m)$.) Hence, $D(\x, \mathcal{X}^n \cup \mathcal{Y}^m) \leq P_{n+m}(\u_\x^\top X \leq \u_\x^\top) = m / (n + m)$. While for $\y$, $P_{n+m}(\u^\top X \leq \u^\top \y) \ge m/(n+m)$ for any $\u \in \mathcal{S}^{d-1}$. Finally, we obtain
\begin{eqnarray*}
    D(\y, \mathcal{X}^n \cup \mathcal{Y}^m) = \sup_{\z \in \textbf{cov}(\mathcal{X}^n \cup \mathcal{Y}^m) \setminus \textbf{cov}(\mathcal{X}^n)} D(\z, \mathcal{X}^n \cup \mathcal{Y}^m).
\end{eqnarray*}
This, together with \eqref{inConX}, implies that $\y \in \mathcal{M}(\mathcal{X}^n \cup \mathcal{Y}^m)$ when $m = n\lambda^*(\mathbf{X}_{\u_0}^n)$.
Note that: (a) $T^*(\mathcal{X}^n \cup \mathcal{Y}^m)$ is by definition the average of all points contained in $\mathcal{M}(\mathcal{X}^n \cup \mathcal{Y}^m)$, (b) $\y$ is arbitrary, it may belong to any bounded region. Hence, $n\lambda^*(\mathbf{X}_{\u_0}^n)$ such $y$ can make $T^*(\mathcal{X}^n \cup \mathcal{Y}^m)$ outside the convex hull of $\mathcal{X}^n$, and in turn break down $T^*$.
\medskip

\emph{This completes the first part of this theorem. Now we proceed to the second part.} By Lemma 5, for $\forall \y \in \mathcal{R}^d \setminus \textbf{cov}(\mathcal{X}^n)$, there must exist a $\u_\y \in \mathcal{S}^{d-1}$ such that $\mathbb{A}_{\u_\y}^\top \y \in \mathcal{M}(\mathbf{X}_{\u_\y}^n)$. Using this and Lemma 6, there $\exists \z$ on the boundary of $\textbf{cov}(\mathcal{X}^n)$, and hence $\z \in \textbf{cov}(\mathcal{X}^n)$, such that
\begin{eqnarray*}
    D(\z, \mathcal{X}^n \cup \mathcal{Y}^m) &\ge& \frac{\min\{n \lambda^*(\mathbf{X}_{\u_\y}^n), ~ m + 1\}}{n + m}\\
    &\ge& \frac{\min\{n \lambda^*(\mathbf{X}_{\u_0}^n), ~ m + 1\}}{n + m}\\
     &>& \frac{m}{n + m} = D(\y, \mathcal{X}^n \cup \mathcal{Y}^m) = \sup_{\z \in \textbf{cov}(\mathcal{X}^n \cup \mathcal{Y}^m) \setminus \textbf{cov}(\mathcal{X}^n)} D(\z, \mathcal{X}^n \cup \mathcal{Y}^m)
\end{eqnarray*}
when $m \leq n \lambda^*(\mathbf{X}_{\u_0}^n) - 1$ for $\u_0$ given in the earlier paragraph of the proof of this theorem. Hence, $T^*(\mathcal{X}^n \cup \mathcal{Y}^m) \in \textbf{cov}(\mathcal{X}^n)$. That is, \emph{less than} $n \lambda^*(\mathbf{X}_{\u_0}^n)$ repetitions of an arbitrary $\y$ \emph{could not} break down $T^*$, no matter where $y$ locates at.
\medskip

This completes the whole proof of this theorem. \hfill $\Box$
\medskip

\textbf{Remark 3.1}. When $d = 2$, $\lambda^*(\mathbf{X}_{\u_0}^n) = \lceil n / 2 \rceil$ for $\u_0 \in \mathcal{S}^{1}$ given in this theorem. Hence, Theorem 1 reduces to the following special case:
\begin{eqnarray*}
    \varepsilon(T^*, \mathcal{X}^n) = \frac{\lceil\frac{n}{2}\rceil}{n + \lceil\frac{n}{2}\rceil}.
\end{eqnarray*}

The key step of Theorem 1 is to locate the new maximizers of Tukey's halfspace depth function after adding $\mathcal{Y}^m$ to $\mathcal{X}^n$. Considering the $\mathbb{A}_\u$-projections of the original observations is a helpful way to achieve this goal of identifying the maximizer. It turns out that the point $\z$ on the boundary of $\textbf{cov}(\mathcal{X}^n)$, that determines a unit vector $\u_\y = (\y - \z) / \|\y - \z\|$ such that the $\mathbb{A}_{\u_\y}$-projections of $\y$ lies in the interior of $\mathcal{M}(\mathbf{X}_{\u_\y}^n)$, plays a key role in the whole proof of Theorem 1.
\medskip

\begin{figure}[H]
\centering
\includegraphics[angle=0,width=5.0in]{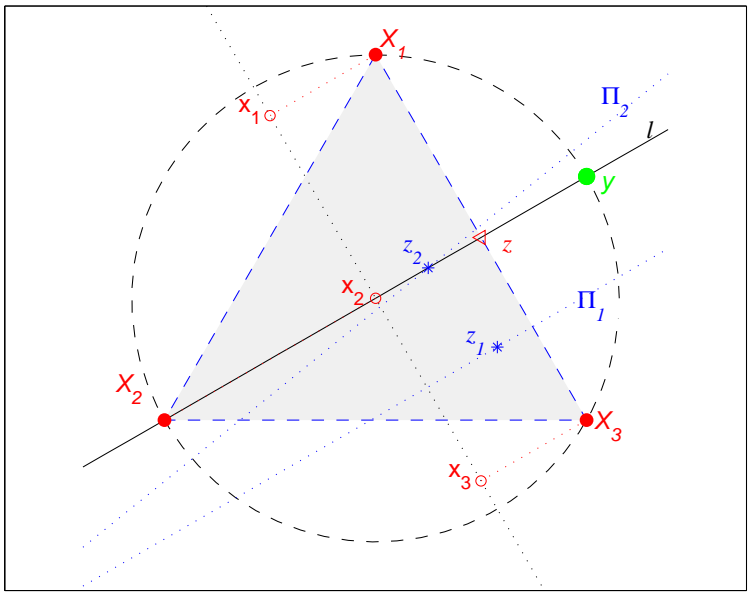}
\caption{Shown is a 2-dimensional illustration for Theorem 1.}
\label{fig:Illustration2D}
\end{figure}

To gain an intuitive understanding of this, we provide a 2-dimensional illustration in Figure \ref{fig:Illustration2D}, where $X_1, X_2, X_3$ denote the data points, and $\mathbf{x}_1, \mathbf{x}_2, \mathbf{x}_3$ the corresponding $\mathbb{A}_{\u_\y}$-projections. When $m=1$, the Tukey depth of the point $\z$ with respect to $\mathcal{X}^n \cup \mathcal{Y}^m = \{X_1, X_2, X_3, \y\}$ is clearly $1/2$, greater than that of any point outside the convex hull of $\{X_1, X_2, X_3\}$. On the other hand, the depth of any $\x \in \textbf{cov}(\mathcal{X}^n) \setminus \{\z\}$ is smaller than $1/2$ (see $\z_1$, $\z_2$ for example).
\medskip

Note that when $\u_0 \in \mathcal{S}^{d-1}$ satisfies \eqref{eqnu0}, $\mathbf{X}_{\u_0}^n$ is in general position from Lemma 1. Hence, relying on Proposition 2.3 in DG92, Theorem 1 in \cite{LLZ2015} and Theorem 1 above, we can easily obtain the following proposition. Since the proof is trivial, we omit it here.
\medskip

\textbf{Proposition 1}. Suppose $\mathcal{X}^n$ is in general position. When $d \ge 2$, the FSBP of Tukey's halfspace median $T^*$ satisfies that
\begin{eqnarray*}
    \frac{\left\lceil \frac{n}{d} \right\rceil}{n + \left\lceil \frac{n}{d} \right\rceil}\leq \varepsilon(T^*, \mathcal{X}^n) \leq
    \left\{ \begin{array}{lcl}
    \frac{\lfloor \frac{n - d + 3}{2}\rfloor}{n + \lfloor \frac{n - d + 3}{2}\rfloor},\  & & \text{if } \exists \u_0 \in \mathcal{S}^{d-1} \text{ satisfying } \eqref{eqnu0},\\
     & &\text{ and } \mathcal{M}(\mathbf{X}_{\u_0}^n) \text{ is singleton},\\
     & &\\
    \frac{\lfloor \frac{n - d + 2}{2}\rfloor}{n + \lfloor \frac{n - d + 2}{2}\rfloor}, \ & &\text{otherwise}.\\
    \end{array}
    \right.
\end{eqnarray*}

\textbf{Remark 3.2}. For $d = 2$, (a) when $n$ is even, $\mathcal{M}(\mathbf{X}_{\u_0}^n)$ is of affine dimension 1, we have $\lfloor \frac{n - d + 2}{2}\rfloor = \lceil\frac{n}{2}\rceil$; (b) when $n$ is odd, $\mathcal{M}(\mathbf{X}_{\u_0}^n)$ is singleton, we have $\lfloor \frac{n - d + 3}{2}\rfloor = \lceil\frac{n}{2}\rceil$. Both scenarios indicate that $\varepsilon(T^*, \mathcal{X}^n)$ attains the upper bound $\lceil\frac{n}{2}\rceil / (n + \lceil\frac{n}{2}\rceil)$.
\medskip

\begin{figure}[H]
\begin{center}
	\subfigure[The scatter plot of the data set.]{
	\includegraphics[angle=0,width=2.5in]{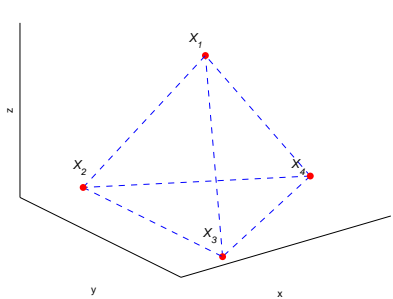}
	\label{fig:DataSet}
	}\quad\quad\quad\quad
	\subfigure[The first scenario of the $\mathbb{A}_u$-projections.]{
	\includegraphics[angle=0,width=2.5in]{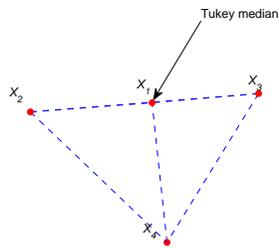}
	\label{fig:Case1}
	}\quad
	\subfigure[The second scenario of the $\mathbb{A}_\u$-projections.]{
	\includegraphics[angle=0,width=2.5in]{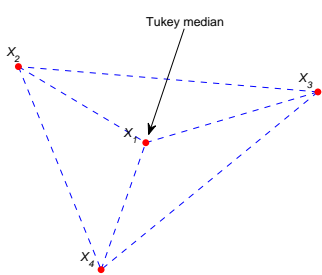}
	\label{fig:Case2}
	}
	\subfigure[The third scenario of the $\mathbb{A}_\u$-projections.]{
	\includegraphics[angle=0,width=2.5in]{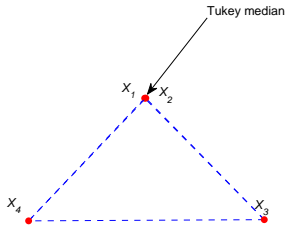}
	\label{fig:Case3}
	}\quad
	\subfigure[The fourth scenario of the $\mathbb{A}_\u$-projections.]{
	\includegraphics[angle=0,width=2.5in]{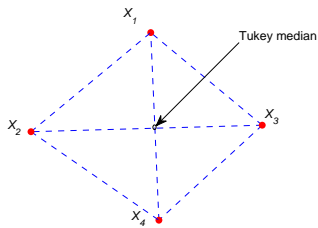}
	\label{fig:Case4}
	}
\caption{Shown is an example for the upper bound for Proposition 1.}
\label{fig:TM3DExam}
\end{center}
\end{figure}

Both the upper and low bound given in Proposition 1 is \emph{attained} if the data set is strategically choosed. Let's first see an illustration for the upper bound. The data points are plotted in \ref{fig:DataSet}. The scatter plot of the $\mathbb{A}_\u$-projections of this data set has four scenarios, though, as shown in Figures \ref{fig:Case1}-\ref{fig:Case4}. The maximum Tukey depth $\lambda^*(\mathbf{X}_\u^n)$ is equal to $1/2$ for any $\u \in \mathcal{S}^{d-1}$, nevertheless. Clearly, $1/2 = \left\lfloor \frac{4 - 3 + 3}{2} \right\rfloor / 4$, and hence $\varepsilon(T^*, \mathcal{X}^n) = 1 / 3$ attains the upper bound given Proposition 1 for this data set.
\medskip

\begin{figure}[H]
\begin{center}
	\subfigure[The scatter plot of the original data set.]{
	\includegraphics[angle=0,width=2.5in]{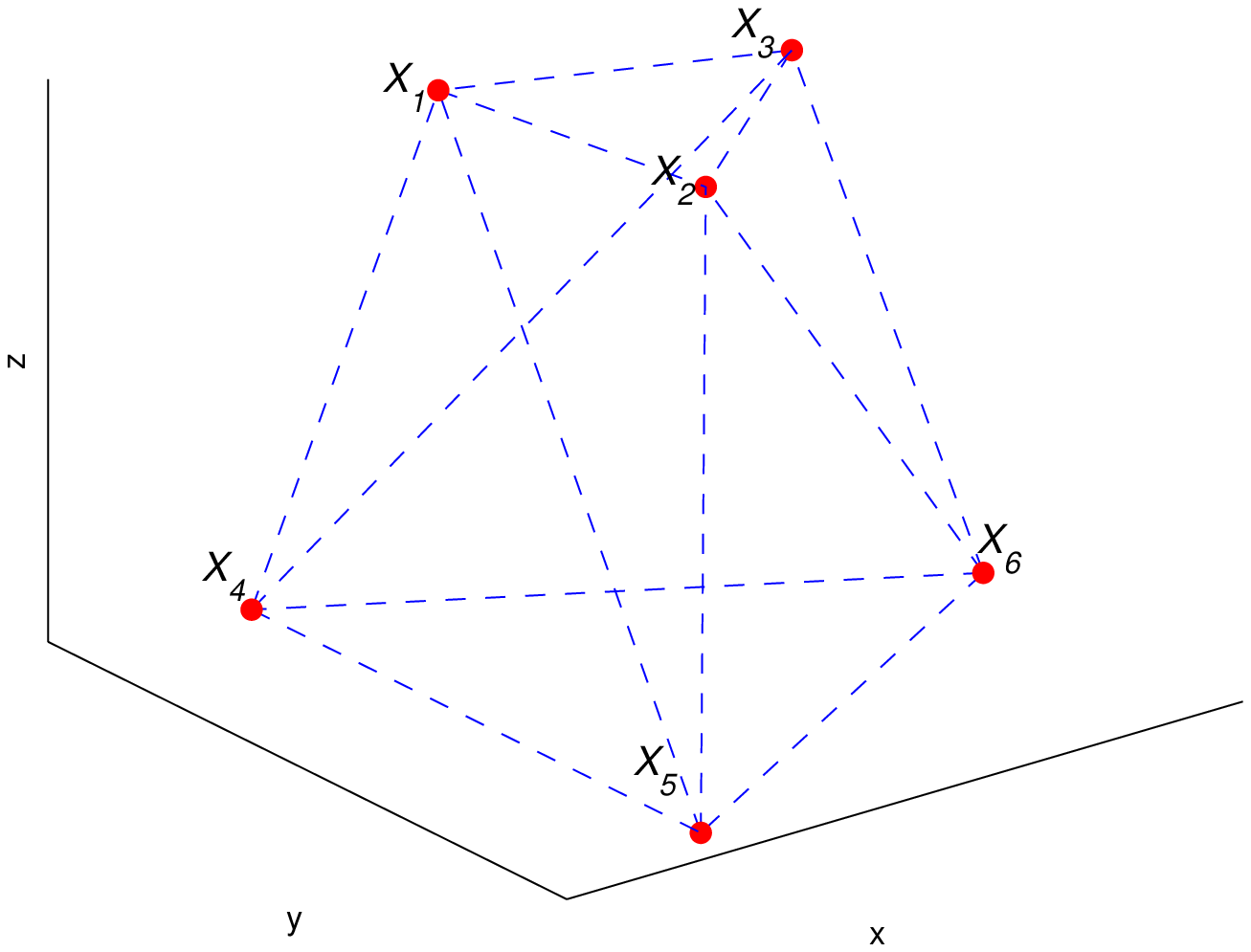}
	\label{fig:SCA1}
	}\quad
	\subfigure[The scatter plot of the $\mathbb{A}_\u$-projections.]{
	\includegraphics[angle=0,width=2.5in]{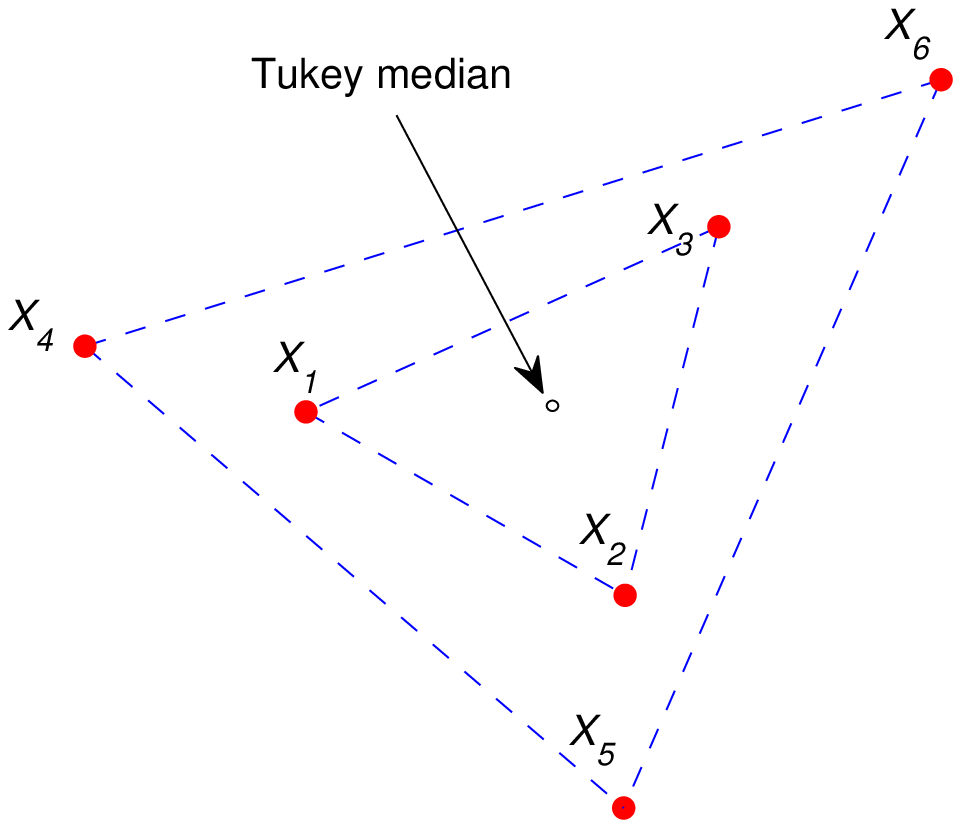}
	\label{fig:PRO1}
	}
\caption{Shown is an example for the low bound for Proposition 1.}
\label{fig:Low}
\end{center}
\end{figure}

As to the low bound, we have an example shown in Figure \ref{fig:SCA1}. Since we can find a $\u$ such that the $\mathbb{A}_\u$-projections of the original data set is a data set of points at the vertices of a collection of nested simplices; See Figure \ref{fig:PRO1}. The maximum Tukey depth with respect to these projections is only $2/6 = 1 / d$ when $d = 3$. Hence, similar to DG92, the low bound of Proposition 1 is also attained, with $\varepsilon(T^*, \mathcal{X}^n) = 1 / 4$ for this example.
\medskip

Compared to the asymptotic result $1/3$, Proposition 1 indicates that the dimension $d$ indeed affects the finite sample breakdown point robustness of Tukey's halfspace median. In detail, when $d$ increases, $\varepsilon(T^*, \mathcal{X}^n)$ tends to decrease for fixed $n$. In fact, the true FSBP of Tukey's halfspace median may be less than $1 / 3$ under the \textbf{IGP} assumption, and this gap may be very great in practice when $d$ is large relative to $n$.

\section{Concluding remarks}
\label{Conclusion}

In the literature, it has long been a open question as to the exact finite sample breakdown point of Tukey's halfspace median. In this paper, we resolved this question through taking account of the $\mathbb{A}_\u$-projections of the original observations when they are in general position. A precise result was provided for \emph{fixed} sample size $n$. The current results revealed that, complimenting the asymptotic result (1/3) obtained by DG92, the finite sample breakdown point robustness of $\HM$ may be affected greatly by the dimension $d$, especially when $d$ is large relative to $n$. Since many offsprings, such as regression depth and multiple output regression, originated directly from Tukey's halfspace depth function with the finite sample breakdown point of their median-like estimators unsolved, we wish that the developed results have the potential to facilitate the investigation of their finite sample breakdown point robustness.
\medskip

Observe that $\inf_{\u \in \mathcal{S}^{d-1}} \lambda^*(\mathbf{X}_\u^n)$ involves an infinite number of maximum Tukey depths $\lambda^*(\mathbf{X}_\u^n)$. It computation is not trivial, and would be very time-consuming. Quite fortunately, there has been much progress in the computation of Tukey's halfspace median and its related depth; See, for example, \cite{RR1998}, \cite{SR2000} and \cite{LLZ2015} and reference therein.

\section*{Acknowledgements}

The research of the first author is supported by National Natural Science Foundation of China (Grant No.11461029, 61263014, 61563018), NSF of Jiangxi Province (No.20142BAB211014, 20143ACB21012,  20132BAB201011, 20151BAB211016), and the Key Science Fund Project of Jiangxi provincial education department (No.GJJ150439, KJLD13033, KJLD14034). Wang's research was supported by the National Science Fund for Distinguished Young Scholars in China (10725106), the National Natural Science Foundation of China (General program 11171331 and Key program 11331011), a grant from the Key Lab of Random Complex Structure and Data Science, CAS and Natural Science Foundation of SZU.

\end{document}